\documentclass{article}

\title{\LARGE \textbf{On large $k$-ended trees in connected graphs}}
\author{Zh.G. Nikoghosyan\footnote{G.G. Nicoghossian (up to 1997)}  }

\begin{document}

\maketitle

\begin{abstract}

A vertex of degree one is called an end-vertex, and an end-vertex of a tree is called a leaf. A tree with at most $k$ leaves is called a $k$-ended tree. For a positive integer $k$, let $t_k$ be the order of a largest $k$-ended tree.   Let $\sigma_m$ be the  minimum degree sum of an independent set of $m$  vertices. The main result (Theorem 2) provides a lower bound for $t_{k+1}$ in terms of $\sigma_m$ and relative orders: if $G$ is a connected graph and $k$, $\lambda$, $m$ are positive integers with $2\le m\le\min\{k,\lambda\}+1$ then either $t_{k+1} \ge \sigma_m +\lambda(k-m+1)+1$ or $t_k\ge t_{k+1}-\lambda+1$.   \\


\end{abstract}

\section{Introduction}

Throughout this article we consider only finite undirected graphs without loops or multiple edges. The set of vertices of a graph $G$ is denoted by $V(G)$ and the set of edges by $E(G)$.  A good reference for any undefined terms is $\cite{[1]}$.

For a graph $G$, we use $n$, $\delta$ and $\alpha$ to denote the order (the number of vertices), the minimum degree and the independence number of $G$, respectively.  If $\alpha\ge k$ for some integer $k$, let $\sigma_k$ be the minimum degree sum of an independent set of $k$ vertices; otherwise we let $\sigma_k=+ \infty$. For a subset $S\subseteq V(G)$, we denote by $G[S]$ the subgraph of $G$ induced by $S$.

If $Q$ is a path or a cycle in a graph $G$, then the order of $Q$, denoted by $|Q|$, is $|V(Q)|$. Each vertex and edge in $G$ can be interpreted as simple cycles of orders 1 and 2, respectively. The graph $G$ is hamiltonian if $G$ contains a Hamilton cycle, i.e. a cycle containing every vertex of $G$. A cycle (path, tree) $Q$ of $G$ is said to be dominating if $V(G-Q)$ is an independent set of vertices.

We write a cycle $Q$ with a given orientation by $\overrightarrow{Q}$. For $x,y\in V(Q)$, we denote by $x\overrightarrow{Q}y$ the subpath of $Q$ in the chosen direction from $x$ to $y$.  For $x\in V(Q)$, we denote the $h$-th successor  and the $h$-th predecessor of $x$ on $\overrightarrow{Q}$ by $x^{+h}$ and $x^{-h}$, respectively. We abbreviate $x^{+1}$ and $x^{-1}$ by $x^+$ and $x^-$, respectively. We say that vertex $z_1$ precedes vertex $z_2$ on $\overrightarrow{Q}$ if $z_1$, $z_2$ occur on $\overrightarrow{Q}$ in this order, and indicate this relationship by $z_1\prec z_2$.

A vertex of degree one is called an end-vertex, and an end-vertex of a tree is usually called a leaf. The set of end-vertices of $G$ is denoted by $End(G)$. A spanning tree is called independence if $End(G)$ is independent in $G$. A branch vertex of a tree is a vertex of degree at least three. The set of branch vertices of a tree $T$ will be denoted by $B(T)$. A tailing of a tree $T$ is a path in $T$ connecting any end-vertex of $T$ to a predecessor of a nearest branch vertex. For a positive integer $k$, a tree $T$ is said to be a $k$-ended tree if $|End(T)|\le k$. A Hamilton path is a spanning 2-ended tree. A Hamilton cycle can be interpreted as a spanning 1-ended tree. In particular, $K_2$ can be interpreted as a hamiltonian graph and as a 1-ended tree. We denote by $t_k$ the order of a largest $k$-ended tree in $G$. In particular, $t_1$ is the order of a longest cycle (the circumference), and $t_2$ is the order of a longest path in $G$.

We first present two simple properties of $k$-ended trees with relative orders $t_k\ge t_{k+1}-\lambda+1$ when $\lambda\in\{1,2\}$. For $\lambda=1$, the following can be checked easily.   \\

\noindent\textbf{Proposition 1}. Let $G$ be a connected graph and $k$ a positive integer. Then $G$ has a spanning $k$-ended tree if and only if $t_k=t_{k+1}$.\\

For $\lambda=2$, we have the dominating version of Proposition 1.  \\

\noindent\textbf{Proposition 2}. Let $G$ be a connected graph with $t_k\ge t_{k+1}-1$ for some positive integer $k$. Then every largest $k$-ended tree in $G$ is a dominating tree.\\

\noindent\textbf{Proof}. Let $G$ be a connected graph with $t_k\ge t_{k+1}-1$ for some $k\ge 1$ and $T_k$ a largest $k$-ended tree in $G$. Suppose the contrary, that is $G-T_k$ contains a component $H$ with $|H|\ge 2$. Now it is easy to construct a $(k+1)$-ended tree $T_{k+1}$ that contains all vertices of $T_k$ and at least 2 vertices of $H$. Then $t_{k+1}\ge |T_{k+1}|\ge t_{k}+2$, contradicting $t_k\ge t_{k+1}-1$.    \ \ \ \   $\triangle$ \\

Our starting point is the earliest degree sum condition for a graph to be hamiltonian due to Ore \cite{[5]}. \\

\noindent\textbf{Theorem A} \cite{[5]}. Every graph with $\sigma_2\ge n$ is hamiltonian.\\

The analog of Theorem A for Hamilton paths follows easily.\\

\noindent\textbf{Theorem B} \cite{[5]}. Every graph with $\sigma_2\ge n-1$ has a Hamilton path.\\

In 1971, Las Vergnas \cite{[3]} gave a degree condition that guarantees that any forest in $G$ of limited size and with a limited number of leaves can be extended to a spanning tree of $G$ with a limited number of leaves in an appropriate sense. This result implies as a corollary a degree sum condition for the existence of a tree with at most $k$ leaves including Theorem A and Theorem B as special cases for $k=1$ and $k=2$, respectively.   \\

\noindent\textbf{Theorem C} \cite{[2]}, \cite{[3]}, \cite{[4]}. Let $G$ be a connected graph with $\sigma_2\ge n-k+1$ for some positive integer $k$. Then $G$ has a spanning $k$-ended tree.\\

However, Theorem C was first openly formulated and proved in 1976 by the author \cite{[4]} and was reproved in 1998 by  Broersma and Tuinstra \cite{[2]}.

In this paper we first present a non-degree sum condition for relative orders $t_k$ and $t_{k+1}$.\\

\noindent\textbf{Theorem 1}. Let $G$ be a connected graph and
let $k$ and $\lambda$ be positive integers with $k\ge 2$.
If $\lambda\ge t_{k+1}/(k+1)$, then $t_k\ge t_{k+1}-\lambda+1$.\\

Since $n\ge t_{k+1}$, Theorem 1 implies the following immediately.\\

\noindent\textbf{Corollary 1}. Let $G$ be a connected graph and
let $k$ and $\lambda$ be positive integers with $k\ge 2$.
If $\lambda\ge n/(k+1)$, then $t_k\ge t_{k+1}-\lambda+1$.\\

The next relation follows from Theorem 1 for a special case when $\lambda=\lfloor t_{k+1}/(k+1)\rfloor$.\\

\noindent\textbf{Corollary 2}. Let $G$ be a connected graph.
Then for each integer $k\ge 2$,
$$
t_k\ge \frac{k}{k+1}t_{k+1}+\frac{1}{k+1}.
$$

The next two results of this paper provide a generalized degree sum conditions for trees with few leaves in connected graphs including Theorems A,\ B,\ C,\ D as special cases. \\

\noindent\textbf{Theorem 2}. Let $G$ be a connected graph and let $k$, $\lambda$, $m$ be positive integers with $2\le m\le\min\{k,\lambda\}+1$. Then either
$$
t_{k+1} \ge \sigma_m +\lambda(k-m+1)+1
$$
or $t_k\ge t_{k+1}-\lambda+1$.  \\

\noindent\textbf{Theorem 3}. Let $G$ be a connected graph and let $k$, $\lambda$, $m$ be positive integers with $m\le\min\{k,\lambda\}+1$. If
$$
\sigma_m\ge t_{k+1}-\lambda(k-m+1)
$$
then $t_k\ge t_{k+1}-\lambda+1$.  \\

Theorem 3 follows from Theorem 2 immediately. The graph
$$
G_1=(k+1)K_{\lambda}+K_1\equiv (mK_{\lambda}\cup (k-m+1)K_{\lambda})+K_1
$$
 shows that the condition $\sigma_m\ge t_{k+1}-\lambda(k-m+1)$ in Theorem 3 cannot be relaxed to $\sigma_m\ge t_{k+1}-\lambda(k-m+1)-1$.

 Next, the graph
$$
G_2=(k+1)K_{\lambda-1}+K_1\equiv (mK_{\lambda-1}\cup (k-m+1)K_{\lambda-1})+K_1
$$
shows that the conclusion $t_k\ge t_{k+1}-\lambda+1$
 in Theorem 3 cannot be strengthened to $t_k\ge t_{k+1}-\lambda+2$
 when $m\le k$. If $m=k+1$ then for this purpose we can use  the
 graph $(k+2)K_{k-1}+K_2$ when $k\ge 2$, and the complete bipartite
 graph $K_{r,r}$ when $k=1$.  Thus, Theorem 3 is best possible.

Theorem 3 implies a number of results in more popular terminology, including Ore-type versions, as well as their spanning $k$-ended and dominating $k$-ended versions.\\

\noindent\textbf{Corollary 3} (Theorem 3, $n\ge t_{k+1}$).

\noindent Let $G$ be a connected graph and let $k$, $\lambda$, $m$ be positive integers with $m\le\min\{k,\lambda\}+1$. If
$$
\sigma_m\ge n-\lambda(k-m+1)
$$
then $t_k\ge t_{k+1}-\lambda+1$.  \\

\noindent\textbf{Corollary 4} (Theorem 3, $m=k+1=\lambda+1$).

\noindent Let $G$ be a connected graph with $\sigma_{k+1}\ge t_{k+1}$ for some positive integer $k$. Then $t_k\ge t_{k+1}-k+1$.\\

\noindent\textbf{Corollary 5} (Theorem 3, $m=2$).

\noindent Let $G$ be a connected graph with $\sigma_2\ge t_{k+1}-\lambda (k-1)$ for some positive integers $\lambda, k$. Then $t_k\ge t_{k+1}-\lambda+1$.\\

\noindent\textbf{Corollary 6} (Theorem 3, $m=2,\ \lambda=1$).

\noindent Let $G$ be a connected graph with $\sigma_2\ge t_{k+1}-k+1$ for some positive integer $k$. Then $t_k\ge t_{k+1}$.\\

\noindent\textbf{Corollary 7} \cite{[2]}, \cite{[3]}, \cite{[4]} (Theorem 3, $m=2,\ \lambda=1$).

\noindent Let $G$ be a connected graph with $\sigma_2\ge n-k+1$ for some positive integer $k$. Then $G$ has a spanning $k$-ended tree.\\

\noindent\textbf{Corollary 8} (Theorem 3, $m=2,\ \lambda=2$).

\noindent Let $G$ be a connected graph with $\sigma_2\ge t_{k+1}-2k+2$ for some positive integer $k$. Then $t_k\ge t_{k+1}-1$.\\

\noindent\textbf{Corollary 9} (Theorem 3, $m=2,\ \lambda=2$).

\noindent Let $G$ be a connected graph with $\sigma_2\ge n-2k+2$ for some positive integer $k$. Then $G$ has a dominating $k$-ended tree.\\

\noindent\textbf{Corollary 10} (Theorem 3, $m=2,\ k=\lambda=1$).

\noindent Let $G$ be a connected graph with $\sigma_{2}\ge t_{2}$. Then $t_1\ge t_{2}$.\\

\noindent\textbf{Corollary 11} \cite{[5]} (Theorem 3, $m=2,\ k=\lambda=1$).

\noindent Let $G$ be a connected graph with $\sigma_{2}\ge n$. Then $G$ is hamiltonian.\\

\noindent\textbf{Corollary 12} (Theorem 3, $m=2,\ \lambda=1,\ k=2$).

\noindent Let $G$ be a connected graph with $\sigma_{2}\ge t_{3}-1$. Then $t_2\ge t_3$.\\

\noindent\textbf{Corollary 13} \cite{[5]} (Theorem 3, $m=2,\ \lambda=1,\ k=2$).

\noindent Let $G$ be a connected graph with $\sigma_{2}\ge n-1$. Then $G$ has a Hamilton path.\\

\noindent\textbf{Corollary 14} (Theorem 3, $m=2,\ k=\lambda=2$).

\noindent Let $G$ be a connected graph with $\sigma_{2}\ge t_{3}-2$. Then $t_2\ge t_3-1$.\\

\noindent\textbf{Corollary 15} (Theorem 3, $m=2,\ k=\lambda=2$).

\noindent Let $G$ be a connected graph with $\sigma_{2}\ge n-2$. Then $G$ has a dominating path.\\

\noindent\textbf{Corollary 16} (Theorem 3, $m=3$).

\noindent Let $G$ be a connected graph with $\sigma_{3}\ge t_{k+1}-\lambda (k-2)$ for some integers $k\ge2$ and $\lambda\ge2$. Then $t_k\ge t_{k+1}-\lambda+1$.\\

\noindent\textbf{Corollary 17} (Theorem 3, $m=3,\ \lambda=2$).

\noindent Let $G$ be a connected graph with $\sigma_{3}\ge t_{k+1}-2k+4$ for some integer $k\ge2$. Then $t_k\ge t_{k+1}-1$.\\

\noindent\textbf{Corollary 18} (Theorem 3, $m=3,\ \lambda=2$).

\noindent Let $G$ be a connected graph with $\sigma_{3}\ge n-2k+4$ for some integer $k\ge2$. Then $G$ has a dominating $k$-ended tree.\\

\noindent\textbf{Corollary 19} (Theorem 3, $m=3,\ k=\lambda=2$).

\noindent Let $G$ be a connected graph with $\sigma_{3}\ge t_{3}$.  Then $t_2\ge t_{3}-1$.\\

\noindent\textbf{Corollary 20} (Theorem 3, $m=3,\ k=\lambda=2$).

\noindent Let $G$ be a connected graph with $\sigma_{3}\ge n$.  Then $G$ has a dominating path.\\

\section{Proofs}

\noindent\textbf{Proof of Theorem 1}. For a connected graph $G$ and positive integers $\lambda$ and $k\ge 2$, let $T_{k+1}$ be a $(k+1)$-ended tree in $G$ and let $A_1,A_2,...,A_{k+1}$ be the tailings of $T_{k+1}$. Clearly, $T_{k+1}-A_i$ is a $k$-ended tree in $G$ for each $i\in \{1,2,...k+1\}$. If $|A_i|\le (t_{k+1}-1)/(k+1)$ for some $i\in\{1,2,...,k+1\}$, then
$$
t_k\ge|T_{k+1}-A_i|=|T_{k+1}|-|A_i|
$$
$$
\ge t_{k+1}-\frac{t_{k+1}-1}{k+1}=t_{k+1}-\lambda+\frac{1}{k+1},
$$
implying that $t_k\ge t_{k+1}-\lambda+1$.

Now let $|A_i|\ge t_{k+1}/ (k+1)$ for each $i\in \{1,2,...k+1\}$. Since $k\ge 2$ and $G$ is connected, $T_{k+1}$ has  a branch vertex $x$. By the definition, $x\not\in A_i$ $(i=1,2,...,k+1)$. Then
$$
t_{k+1}\ge \sum_{i=1}^{k+1}|A_i|+|\{x\}|\ge t_{k+1}+1,
$$
a contradiction.          \qquad \rule{7pt}{6pt}    \\

\noindent\textbf{Proof of Theorem 2}. Let $G$ be a connected graph and let $k$, $\lambda$, $m$ be positive integers with
$2\le m\le \min\{k,\lambda\}+1$. If $t_k\ge t_{k+1}-\lambda+1$ then we are done. Let
$$
t_k\le t_{k+1}-\lambda.                  \eqno{(1)}
$$
We shall prove that
$$
t_{k+1}\ge \sigma_m+\lambda(k-m+1)+1.               \eqno{(2)}
$$
Let $T_{k+1}$ be a $(k+1)$-ended tree in $G$ and $T_m$ be an $m$-ended subtree of $T_{k+1}$. Assume that  \\

(i) \ $T_{k+1}$ is chosen so that $|E(T_{k+1})|$ is as large as possible,\\

(ii) \ $T_{k+1}$ is chosen so that $|E(T_{m})|$ is as large as possible, subject to (i).\\

By the definition, $|T_{k+1}|=t_{k+1}$.\\

\noindent\textbf{Claim 1}. $|End(T_{k+1})|= k+1\ge2$.\\

\noindent\textbf{Proof}. Assume the contrary, that is $|End(T_{k+1})|\le k$, implying that $T_{k+1}$ is a $k$-ended tree. Since $\lambda\ge 1$, we have
$$
t_k\ge |T_{k+1}|=t_{k+1}\ge t_{k+1}-\lambda+1,
$$
contradicting (1). Hence,  $|End(T_{k+1})|= k+1$.  Recalling also that $k\ge 1$, we have $|End(T_{k+1})|\ge 2$.    \ \ \ \   $\triangle$ \\

\noindent\textbf{Claim 2}. $T_{k+1}$ is an independence tree.\\

\noindent\textbf{Proof}. If two of the end-vertices of $T_{k+1}$ are joined by an edge $e$, then $T_{k+1}+e$ has a unique cycle $C$. If $C$ is a Hamilton cycle, then $T_{k+1}$ is a 1-ended tree, contradicting Claim 1. Otherwise at least one vertex $v$ of $C$ has a degree at least three in $T_{k+1}+e$. Deleting one of the edges of $C$ incident with $v$ results in a $k$-ended tree $T_k$ of order $|T_{k+1}|$. Then
$$
t_k\ge |T_k|=|T_{k+1}|=t_{k+1}\ge t_{k+1}-\lambda+1,
$$
contradicting (1). Hence, $T_{k+1}$ is an independence tree.   \ \ \ \   $\triangle$ \\

\noindent\textbf{Claim 3}. If $L$ is a tailing of a $(k+1)$-ended tree $T$ in $G$ with $|T|=t_{k+1}$, then $|L|\ge \lambda$.\\

\noindent\textbf{Proof}.  Assume the contrary, that is $|L|\le \lambda-1$ for some tailing $L$ of $T$. Since $T-L$ is a $k$-ended tree, we have
$$
t_k\ge |T-L|= |T|-|L|\ge t_{k+1}-\lambda+1,
$$
contradicting (1).    \ \ \ \   $\triangle$ \\

\noindent\textbf{Claim 4}. If $T$ is a $k$-ended tree in $G$ then $|T|<t_{k+1}$.\\

\noindent\textbf{Proof}. Assume the contrary, that is $|T|\ge t_{k+1}$. Then
$$
t_k\ge|T|\ge t_{k+1}\ge t_{k+1}-\lambda+1,
$$
contradicting (1).    \ \ \ \   $\triangle$ \\

\textbf{Case 1}. $|End(T_{k+1})|=2$.

By Claim 1, $k=1$ and $m=2$, implying that $T_2$ is a longest path in $G$. Put $T_2=v_1v_2...v_f$. By Claim 2, $v_1v_f\not\in E(G)$.  By (i), $N(v_1)\cup N(v_f)\subseteq V(T_2)$. If $d(v_1)+d(v_f)\ge t_2$ then by standard arguments, $G[V(T_2)]$ is hamiltonian, that is  $G[V(T_2)]$ contains a 1-ended tree (cycle) $T_1$ with $|T_1|=|T_2|$, contradicting Claim 4. Otherwise
$$
t_{k+1}=t_2\ge d(v_1)+d(v_f)+1\ge \sigma_2+1=\sigma_m+\lambda(k-m+1)+1.
$$

\textbf{Case 2}. $|End(T_{k+1})|\ge 3$.

Put $End(T_{k+1})=\{\xi_1,\xi_2,...,\xi_{k+1}\}$. By (ii), $End(T_m)\subseteq End(T_{k+1})$. Assume w.l.o.g. that $End(T_{m})=\{\xi_1,\xi_2,...,\xi_{m}\}$.  If $\cup_{i=1}^mN(\xi_i)\not\subseteq V(T_{m})$, then clearly $G$ contains an $m$-ended subtree $T_m^\prime$ with $|E(T_m^\prime)|>|E(T_m)|$, contradicting (ii). Hence,
$$
\bigcup_{i=1}^mN(\xi_i)\subseteq V(T_m).
$$
For each $i\in \{1,...,k+1\}$, let $\overrightarrow{Q}_i=\xi_i\overrightarrow{Q}_iw_i$ be the tailing of $T_{k+1}$ connecting $\xi_i$ to the predecessor $w_i$ of the nearest branch vertex $w_i^\ast$ of $T_{k+1}$. By Claim 3, $|V(Q_i)|\ge \lambda$. Let $w^\prime_i$ be the vertex on $Q_i$ with $|V(\xi_i\overrightarrow{Q}_iw^\prime_i)|= \lambda$. Put
$$
A_i=V(Q_i), \  A_i^\prime=V(\xi_i\overrightarrow{Q}_iw^\prime_i) \ \  (i=1,...,k+1).
$$

\noindent\textbf{Claim 5}. If $|T_m|\ge \sigma_m+1$ then (2) holds.

\noindent\textbf{Proof}. Since $|A_i|\ge |V(Q_i)|\ge\lambda$ for each $i\in\{1,...,k+1\}$, we have
$$
t_{k+1}=|T_{k+1}|=|T_m|+|T_{k+1}-T_m|
$$
$$
\ge \sigma_m+1+\sum_{i=m+1}^{k+1}|A_i|\ge \sigma_m+\lambda (k-m+1)+1,
$$
and (2) holds.    \ \ \ \   $\triangle$ \\

To prove that $|T_m|\ge \sigma_m+1$, which by Claim 5 implies (2), we use mathematical induction on $m$. Assume that $m=2$ (induction basis).
By (ii), $N(\xi_1)\cup N(\xi_2)\subseteq V(T_2)$. If $d(\xi_1)+d(\xi_2)\ge |T_2|$ then by standard arguments, $G[V(T_2)]$ is hamiltonian and we can form a $k$-ended tree $T_{k+1}^\prime$ of order $|T_{k+1}|$, contradicting Claim 4. Otherwise $|T_2|\ge d(\xi_1)+d(\xi_2)\ge \sigma_2+1$. Now suppose that (2) holds for $m-1$, where $m\ge3$.    \\

\noindent\textbf{Claim 6}. Let $\mu\in A_i$ for some $i\in \{1,2,...,k+1\}$. If $\xi_j\mu\in E(G)$ for some $j\in \{1,2,...,k+1\}-\{i\}$, then $|\xi_i\overrightarrow{Q_i}\mu^-|\ge \lambda$ and $|\mu^+\overrightarrow{Q_i}w_i|\ge \lambda$.\\

\noindent\textbf{Proof}. Put
$$
T_{k+1}^\prime=T_{k+1}+\xi_j\mu-w_iw_i^\ast.
$$
By Claim 2, $\mu\not=\xi_i$. Next, we have $\mu\not=w_i$ since otherwise $T_{k+1}^\prime$ is a $k$-ended tree of order $|T_{k+1}|$, contradicting Claim 4. Then $T_{k+1}^\prime$ is a $(k+1)$-ended tree with tailings $\xi_i\overrightarrow{Q_i}\mu^-$  and $\mu^+\overrightarrow{Q_i}w_i$. By Claim 3, $|\xi_i\overrightarrow{Q_i}\mu^-|\ge \lambda$ and $|\mu^+\overrightarrow{Q_i}w_i|\ge \lambda$.    \ \ \ \   $\triangle$ \\

\noindent\textbf{Claim 7}. Let $\mu_1,\mu_2\in A_i$ for some $i\in \{1,2,...,k+1\}$ and let $\mu_1\prec \mu_2$. If $\xi_i\mu_2,\xi_j\mu_1\in E(G)$ for some $j\in \{1,2,...,k+1\}-\{i\}$, then $|\mu_1^+\overrightarrow{Q_i}\mu_2^-|\ge \lambda$.\\

\noindent\textbf{Proof}. Put
$$
T_{k+1}^\prime=T_{k+1}+\xi_i\mu_2+\xi_j\mu_1-\mu_1\mu_1^+-w_iw_i^\ast.
$$
If $\mu_1^+=\mu_2$ then $T_{k+1}^\prime$ is a $k$-ended tree of order $|T_{k+1}|$, contradicting Claim 4. Otherwise $T_{k+1}^\prime$ is a $(k+1)$-ended tree with tailing $\mu_1^+\overrightarrow{Q_i}\mu_2^-$. By Claim 3, $|\mu_1^+\overrightarrow{Q_i}\mu_2^-|\ge \lambda$.    \ \ \ \   $\triangle$ \\

\noindent\textbf{Claim 8}. Let $\mu_1,\mu_2\in A_i$ for some $i\in \{1,2,...,k+1\}$ and let $\mu_1\prec \mu_2$. If $\xi_j\mu_1,\xi_t\mu_2\in E(G)$ for some distinct $j,t\in \{1,2,...,k+1\}-\{i\}$, then $|\mu_1^+\overrightarrow{Q_i}\mu_2^-|\ge \lambda$.\\

\noindent\textbf{Proof}. Put
$$
T_{k+1}^\prime=T_{k+1}+\xi_j\mu_1+\xi_t\mu_2-\mu_1\mu_1^+-w_iw_i^\ast.
$$
If $\mu_1^+=\mu_2$ then $T_{k+1}^\prime$ is a $k$-ended tree of order $|T_{k+1}|$, contradicting Claim 4. Otherwise $T_{k+1}^\prime$ is a $(k+1)$-ended tree with tailing $\mu_1^+\overrightarrow{Q_i}\mu_2^-$. By Claim 3, $|\mu_1^+\overrightarrow{Q_i}\mu_2^-|\ge \lambda$.   \ \ \ \   $\triangle$ \\

\noindent\textbf{Claim 9}. Let $i,j\in\{1,...,m\}$ and $i\not=j$. Then
$$
N_{A_1}^{+(i-1)}(\xi_i)\cap N_{A_1}^{+(j-1)}(\xi_j)=\emptyset.
$$

\noindent\textbf{Proof}. Assume the contrary and let $\mu\in N_{A_1}^{+(i-1)}(\xi_i)\cap N_{A_1}^{+(j-1)}(\xi_j)$.

Assume first that $i\ge 2$ and $j\ge 2$. It follows that $\mu_1^{+(i-1)}=\mu_2^{+(j-1)}=\mu$ for some $\mu_1\in N_{A_1}(\xi_i)$ and $\mu_2\in N_{A_1}(\xi_j)$. Assume w.l.o.g. that $\mu_1<\mu_2$, that is $j<i$. Then
$$
|\mu_1^+\overrightarrow{Q}_1\mu_2^-|=i-j-1\le m-2\le \lambda-1,
$$
contradicting Claim 8.

Now assume that either $i=1$ or $j=1$, say $j=1$. By the hypothesis, $\mu\in N_{A_1}(\xi_1)\cap N_{A_1}^{+(i-1)}(\xi_i)$. It follows that $\mu_1^{+(i-1)}=\mu$ for some $\mu_1\in N_{A_1}(\xi_i)$. Then
$$
|\mu_1^+\overrightarrow{Q}_1\mu^-|=i-2\le m-2\le \lambda-1,
$$
contradicting Claim 7. \ \ \ \   $\triangle$ \\

Since $m\le \lambda+1$, by Claim 6, $N_{A_1}^{+(i-1)}(\xi_i)\subseteq A_1$
 for each $i\in \{1,2,...,m\}$. Next, it is easy to see that $\xi_1\not\in N_{A_1}^{+(i-1)}(\xi_i)$ for each $i\in\{1,...,m\}$. Then by Claim 9,
$$
|A_1|\ge \sum_{i=1}^m|N_{A_1}(\xi_i)|+|\{\xi_1\}|=\sum_{i=1}^m|N_{A_1}(\xi_i)|+1.
$$
By a similar argument, for each $j\in \{1,2,...,m\}$,
$$
|A_j|\ge \sum_{i=1}^m|N_{A_j}(\xi_i)|+1.
$$

Put
$$
A=\bigcup_{i=1}^mA_i.
$$
Clearly,
$$
|A|=\sum_{i=1}^m|A_i|\ge\sum_{i=1}^m|N_A(\xi_i)|+m.               \eqno{(3)}
$$

Let $\Gamma$ be the set of all paths in $T_{m}$ with \\

$(\ast)$ \ \ $M\in \Gamma$ if and only if $|M|\ge 2$ and $V(M)\cap B(T_{k+1})=End(M)$.\\

Let $M_1,M_2,...,M_\pi$ be the elements of $\Gamma$. For each $i\in \{1,2,...,\pi\}$, put
$$
\overrightarrow{M_i}=x_i\overrightarrow{M_i}y_i, \ \  D_i=V(M_i)-\{x_i,y_i\}.
$$
For each $i\in\{1,...,\pi\}$, $T_{m}-D_i$ consists of two connected components, denoted by $T_{m}(x_i)$ and $T_{m}(y_i)$. \\

\textbf{Case 2.1}. $|D_i|\ge \sum_{j=1}^m|N_{D_i}(\xi_j)|+\lambda$  \  \  $(i=1,...,\pi)$.

Put
$$
D=\bigcup_{i=1}^{\pi} D_i, \ \  B^\prime=V(T_m)\cap B(T_{k+1}).
$$
Clearly,
$$
|B^\prime|=\pi+1, \ \ A\cap D=A\cap B^\prime=D\cap B^\prime=\emptyset, \ \ |T_{m}|=|A|+|D|+|B^\prime|.
$$
Since $m\ge3$, we have $|B^\prime|\not=\emptyset$, that is $\pi\ge0$. By the hypothesis,
$$
|D|=\sum_{i=1}^{\pi}|D_i|\ge \sum_{i=1}^m|N_D(\xi_i)|+\pi \lambda.
$$
By (3),
$$
|T_m|=|A|+|D|+|B^\prime|
$$
$$
\ge\left(\sum_{i=1}^m|N_A(\xi_i)|+m\right)+\left(\sum_{i=1}^{m}|N_D(\xi_i)|+\pi\lambda\right)+\pi+1
$$
$$
\ge\left(\sum_{i=1}^m|N_A(\xi_i)|+\sum_{i=1}^{m}|N_D(\xi_i)|+\sum_{i=1}^m|N_{B^\prime}(\xi_i)|\right)-\sum_{i=1}^m|N_{B^\prime}(\xi_i)|
+m+\pi(\lambda+1)+1
$$
$$
\ge\sum_{i=1}^m|N_{T_m}(\xi_i)|- \sum_{i=1}^{m}|N_{B^\prime}(\xi_i)|+m+\pi(\lambda+1)+1.
$$
Observing that
$$
\sum_{i=1}^m|N_{T_m}(\xi_i)|=\sum_{i=1}^md(\xi_i)\ge\sigma_m,
$$
$$
\sum_{i=1}^m|N_{B^\prime}(\xi_i)|\le m|B^\prime|=m(\pi +1),
$$
we get
$$
|T_m|\ge \sigma_m+m+\pi(\lambda+1)-m(\pi+1)+1=\sigma_m+\pi(\lambda-m+1)+1.
$$
Since $\pi\ge0$ and $\lambda\ge m+1$ (by the hypothesis), we have $|T_m|\ge \sigma_m+1$ and (2) holds by Claim 5.  \\

\textbf{Case 2.2}. $|D_i|\le \sum_{j=1}^m|N_{D_i}(\xi_j)|+\lambda-1$ for some  $i\in\{1,...,\pi\}$.

Assume w.l.o.g. that $i=1$, that is
$$
|D_1|\le \sum_{i=1}^m|N_{D_1}(\xi_i)|+\lambda-1.         \eqno{(4)}
$$
Put $T_m(x_1)=H$ and $T_m(y_1)=F$. Assume w.l.o.g. that
$$
End(H)=\{\xi_1,\xi_2,...,\xi_r\}, \ \  End(F)=\{\xi_{r+1},\xi_{r+2},...,\xi_m\},
$$
where $1\le r\le m-1$.\\

\noindent\textbf{Claim 10}. Let $\xi_i\mu\in E(G)$ for some $i\in \{1,...,r\}$, say $i=1$, and $\mu\in D_1\cup V(F)$. If $\mu\in D_1$ then $|x_1^+\overrightarrow{M_1}\mu^-|\ge\lambda$. If $\mu\in V(F)$ then $|D_1|\ge \lambda$.   \\

\noindent\textbf{Proof}. Put
$$
T_{k+1}^\prime=T_{k+1}+\xi_1\mu-x_1x_1^+.
$$
Assume first that $\mu\in D_1$. If $\mu=x_1^+$ then $T_{k+1}^\prime$ is a $k$-ended tree of order $|T_{k+1}|$, contradicting Claim 4. Otherwise $T_{k+1}^\prime$ is a $(k+1)$-ended tree with tailing $x_1^+\overrightarrow{M_1}\mu^-$. By Claim 3, $|x_1^+\overrightarrow{M_1}\mu^-|\ge\lambda$.

 Now let $\mu\in V(F)$. If $x_1^+=y_1$ then $T_{k+1}^\prime$ is a $k$-ended tree of order $|T_{k+1}|$, contradicting Claim 4. Otherwise $T_{k+1}^\prime$ is a $(k+1)$-ended tree with tailing $x_1^+\overrightarrow{M_1}y_1^-$. By Claim 3, $|D_1|\ge \lambda$.  \ \ \ \   $\triangle$ \\

\noindent\textbf{Claim 11}. Let $\mu_1\prec\mu_2$ for some $\mu_1,\mu_2\in D_1$. If $\xi_i\mu_2, \xi_j\mu_1\in E(G)$ for some $i\in \{1,...,r\}$ and $j\in \{r+1,r+2,...,m\}$, then
$$
|\mu_1^+\overrightarrow{M_1}\mu_2^-|\ge\lambda, \ \ |x_1^+\overrightarrow{M_1}\mu_1^-|\ge\lambda, \ \  |\mu_2^+\overrightarrow{M_1}y_1^-|\ge\lambda.
$$

\noindent\textbf{Proof}. Put
$$
T_{k+1}^\prime=T_{k+1}+\xi_i\mu_2+\xi_j\mu_1-x_1x_1^+-\mu_1\mu_1^+.
$$
If $\mu_1^+=\mu_2$ then $T_{k+1}^\prime$ is a $k$-ended tree of order $|T_{k+1}|$, contradicting Claim 4. Otherwise $T_{k+1}^\prime$ is a $(k+1)$-ended tree with tailing $\mu_1^+\overrightarrow{M_1}\mu_2^-$. By Claim 3, $|\mu_1^+\overrightarrow{M_1}\mu_2^-|\ge\lambda$.

Now put
$$
T_{k+1}^{\prime\prime}=T_{k+1}+\xi_i\mu_2+\xi_j\mu_1-x_1x_1^+-y_1y_1^-.
$$
If $x_1^+=\mu_1$ then $T_{k+1}^{\prime\prime}$ is a $k$-ended tree of order $|T_{k+1}|$, contradicting Claim 4. Otherwise $T_{k+1}^{\prime\prime}$ is a $(k+1)$-ended tree with tailing $x_1^+\overrightarrow{M_1}\mu_1^-$. By Claim 3, $|x_1^+\overrightarrow{M_1}\mu_1^-|\ge\lambda$. By a symmetric argument, $|\mu_2^+\overrightarrow{M_1}y_1^-|\ge\lambda$.    \ \ \ \   $\triangle$ \\

\noindent\textbf{Claim 12}. Let $\mu_1\prec\mu_2$ for some $\mu_1,\mu_2\in D_1$ and let $\xi_i\mu_2, \xi_j\mu_1\in E(G)$ for some distinct $i,j\in\{1,...,m\}$. If either $i,j\in\{1,...,r\}$ or $i,j\in\{r+1,r+2,...,m\}$, then $|\mu_1^+\overrightarrow{M_1}\mu_2^-|\ge\lambda$.\\

\noindent\textbf{Proof}. Assume w.l.o.g. that $i,j\in\{1,...,r\}$. Put
$$
T_{k+1}^\prime=T_{k+1}+\xi_i\mu_2+\xi_j\mu_1-x_1x_1^+-\mu_1\mu_1^+.
$$
If $\mu_1^+=\mu_2$ then $T_{k+1}^\prime$ is a $k$-ended tree of order $|T_{k+1}|$, contradicting Claim 4. Otherwise $T_{k+1}^\prime$ is a $(k+1)$-ended tree with tailing $\mu_1^+\overrightarrow{M_1}\mu_2^-$. By Claim 3, $|\mu_1^+\overrightarrow{M_1}\mu_2^-|\ge\lambda$.   \ \ \ \   $\triangle$ \\

\noindent\textbf{Claim 13}. Let $i,j\in \{1,...,m\}$ and $i\not= j$. Then
$$
N_{D_1}^{-(r-i)}(\xi_i)\cap N_{D_1}^{-(r-j)}(\xi_j)=\emptyset.
$$
\noindent\textbf{Proof}. Assume the contrary and let $\mu\in N_{D_1}^{-(r-i)}(\xi_i)\cap N_{D_1}^{-(r-j)}(\xi_j)$.

Assume first that $i\le r$ and $j\ge r+1$.
It follows that $\mu_1^{-(r-j)}=\mu_2^{-(r-i)}=\mu$ for some
$\mu_1\in N_{D_1}(\xi_j)$ and $\mu_2\in N_{D_1}(\xi_i)$.
Since $j>i$, we have $\mu_1\prec\mu_2$. Then
$$
|\mu_1^+\overrightarrow{M_1}\mu_2^-|=(j-r)+(r-i)-1=j-i-1\le m-2\le\lambda-1,
$$
contradicting Claim 11.

Now assume that either $i,j\le r$ or $i,j\ge r+1$, say $i,j\le r$. Assume w.l.o.g. that $i<j$. It follows that $\mu_1^{-(r-j)}=\mu_2^{-(r-i)}=\mu$ for some $\mu_1\in N_{D_1}(\xi_j)$ and $\mu_2\in N_{D_1}(\xi_i)$. Since $i<j$, we have $\mu_1\prec\mu_2$. Then
$$
|\mu_1^+\overrightarrow{M_1}\mu_2^-|=(r-i)-(r-j)-1=j-i-1\le r-2\le m-3\le \lambda-2,
$$
contradicting Claim 12.          \ \ \ \   $\triangle$ \\

\noindent\textbf{Claim 14}. $N_{D_1}^{-(r-i)}(\xi_i)\subseteq D_1$ \ \ $(i=1,...,m)$.

\noindent\textbf{Proof}. If $N_{D_1}^{-(r-i)}(\xi_i)=\emptyset$ then we are done. Let $\mu\in N_{D_1}^{-(r-i)}(\xi_i)$.

Assume first that $i\le r$. It follows that $\mu_1^{-(r-i)}=\mu$ for some $\mu_1\in N_{D_1}(\xi_i)$. By Claim 10, $|x_1^+\overrightarrow{M_1}\mu_1^-|\ge \lambda$. Observing also that $r-i\le r-1\le m-2\le\lambda-1$, we conclude that $\mu\in D_1$.

Now assume that $i\ge r+1$. It follows that $\mu_1^{+(i-r)}=\mu$ for some $\mu_1\in N_{D_1}(\xi_i)$. By Claim 10, $|\mu_1^+\overrightarrow{M_1}y_1^-|\ge \lambda$. On the other hand, $i-r\le m-r\le m-1\le\lambda$. Hence, $\mu\in D_1$.   \ \ \ \   $\triangle$ \\

\textbf{Case 2.2.1}. $\xi_i\mu_1,\xi_j\mu_2\in E(G)$ for some $i\in\{1,...,r\}$, $j\in\{r+1,r+2,...,m\}$, say $i=1$, $j=m$, and $\mu_1,\mu_2\in D_1$.

By Claims 11 and 12, $|D_1|\ge 2\lambda+1$. Put
$$
X_1=\{x_1^{+1},x_1^{+2},...,x_1^{+(\lambda-r+1)}\}, \ \ Y_1=\{y_1^{-1},y_1^{-2},...,y_1^{-(\lambda-m+r)}\}.
$$
Since $|X_1|+|Y_1|=2\lambda-m+1\le2\lambda-2$, we have $X_1\cup Y_1\subseteq D_1$ and $X_1\cap Y_1=\emptyset$.\\

\noindent\textbf{Claim 15}. For each $i\in\{1,...,m\}$,
$$
(X_1\cup Y_1)\cap N_{D_1}^{-(r-i)}(\xi_i)=\emptyset  \  \  (i=1,...,m).
$$
\noindent\textbf{Proof}. Assume the contrary and let $\mu\in (X_1\cup Y_1)\cap N_{D_1}^{-(r-i)}(\xi_i)$ for some $i\in\{1,...,m\}$. Assume w.l.o.g. that $i\le r$. It follows that $\mu_1^{-(r-i)}=\mu$ for some $\mu_1\in N_{D_1}(\xi_i)$. If $\mu\in X_1$ then
$$
|x_1^+\overrightarrow{M_1}\mu_1^-|\le (\lambda-r+1)+(r-i)-1=\lambda-i\le\lambda-1,
$$
contradicting Claim 10. Now let $\mu\in Y_1$. By Claim 11 and by the hypothesis (Case 2.2.1), $|\mu_1^+\overrightarrow{M_1}y_1^-|\ge \lambda$. Then
$$
|Y_1|=\lambda-m+r\ge |\mu\overrightarrow{M_1}y_1^-|=|\mu\overrightarrow{M_1}\mu_1^-|+|\mu_1^+\overrightarrow{M_1}y_1^-|+1\ge r-i+\lambda+1,
$$
implying that $i\ge m+1$, a contradiction.   \ \ \ \   $\triangle$ \\

By Claims 14,15 and 16,
$$
|D_1|\ge \sum_{i=1}^m|N_{D_1}^{-(r-i)}(\xi_i)|+|X_1|+|Y_1|
$$
$$
\ge \sum_{i=1}^m|N_{D_1}(\xi_i)|+2\lambda-m+1\ge \sum_{i=1}^m|N_{D_1}(\xi_i)|+\lambda,
$$
contradicting (4).\\

\textbf{Case 2.2.2}. $\xi_i\mu\in E(G)$ for some $i\in \{1,...,r\}$, say $i=1$, and  $N_{D_1}(\xi_j)=\emptyset$ for each $j\in \{r+1,r+2,...,m\}$.

Assume that $\xi_1$ and $\mu$ are chosen such that $x_1^+\overrightarrow{M_1}\mu$ is as long as possible.\\

\textbf{Case 2.2.2.1}. $|\mu^+\overrightarrow{M_1}y_1^-|\ge\lambda$.

Put
$$
X_1=\{x_1^{+1},x_1^{+2},...,x_1^{+(\lambda-r+1)}\}, \ \ Y_1=\{y_1^{-1},y_1^{-2},...,y_1^{-\lambda}\}.
$$
Clearly, $X_1\cap Y_1=\emptyset$ and $|X_1|+|Y_1|=2\lambda-r+1\ge2\lambda-m+2\ge\lambda+1$. Since $N_{D_1}(\xi_i)=\emptyset$ for each $i\in \{r+1,r+2,...,m\}$, we have $Y_1\cap N_{D_1}^{-(r-i)}(\xi_i)=\emptyset$ for each $i\in\{1,...,m\}$. As in Claim 15, we have also $X_1\cap N_{D_1}^{-(r-i)}(\xi_i)=\emptyset$ for each $i\in\{1,...,r\}$. Recalling also that $N_{D_1}(\xi_i)=\emptyset$ $(i=1,...,m)$, we conclude that for each $i\in\{1,...,m\}$,
$$
(X_1 \cup Y_1)\cap N_{D_1}^{-(r-i)}(\xi_i)=\emptyset.
$$
Then
$$
|D_1|\ge \sum_{i=1}^m|N_{D_1}^{-(r-i)}(\xi_i)|+|X_1|+|Y_1|\ge \sum_{i=1}^m|N_{D_1}(\xi_i)|+\lambda+1,
$$
contradicting (4).\\

\textbf{Case 2.2.2.2}. $|\mu^+\overrightarrow{M_1}y_1^-|\le\lambda-1$.\\

\noindent\textbf{Claim 16}. $N(\xi_i)\subseteq V(H\cup x_1\overrightarrow{M_1}\mu)$ $(i=2,3,...,r)$.

\noindent\textbf{Proof}. Assume the contrary, that is $\xi_i\mu^\prime\in E(G)$ for some $i\in\{2,...,r\}$, say $i=2$, and $\mu^\prime\in V(F\cup \mu^+\overrightarrow{M_1}y_1)$. If $\mu^\prime\in V(\mu^+\overrightarrow{M_1}y_1)$ then by Claim 12, $|\mu^+\overrightarrow{M_1}y_1|\ge \lambda$, contradicting the hypothesis. Let $\mu^\prime\in V(F)$. Put
$$
T_{k+1}^\prime=T_{k+1}+\xi_1\mu-x_1x_1^+.
$$
Since $\xi_2\mu^\prime\in E(G)$, by Claim 10, $|\mu^+\overrightarrow{M_1}y_1|\ge \lambda$, contradicting the hypothesis.  \ \ \ \   $\triangle$ \\

\noindent\textbf{Claim 17}. $N(\xi_i)\subseteq V(F)$ \ $(i=r+1,r+2,...,m)$.

\noindent\textbf{Proof}. Assume the contrary, that is $\xi_i\mu^\prime\in E(G)$ for some $i\in \{r+1,r+2,...,m\}$, say $i=m$, and $\mu^\prime\in V(H)\cup D_1$. By the hypothesis (Case 2.2.2), $\mu^\prime\not\in D_1$. Let $\mu^\prime\in V(H)$. Put
$$
T_{k+1}^\prime=T_{k+1}+\xi_1\mu-x_1x_1^+.
$$
Since $\xi_m\mu^\prime\in E(G)$, by Claim 10, $|\mu^+\overrightarrow{M_1}y_1^-|\ge\lambda$, contradicting the hypothesis.    \ \ \ \   $\triangle$ \\

We steel cannot use induction hypothesis with respect to $H$ or $F$, since possibly $End(H)\cup End(F)\not\subseteq End(T_{k+1})$ and possibly $N(\xi_i)\not\subseteq V(H)$ for some $i\le r$ or $N(\xi_i)\not\subseteq V(F)$ for some $i\ge r+1$. For this purpose, we shall reform $H$ and $F$, as well as $T_{k+1}$ to appropriate $H^\prime$, $F^\prime$ and $T_{k+1}^\prime$ as follows.

If $d_{H}(x_1)\ge 2$ then $H^\prime=H\cup x_1^+\overrightarrow{M_1}\mu\xi_1$ and  $T_{k+1}^\prime=T_{k+1}+\xi_1\mu-x_1x_1^+$. Clearly, $End(H^\prime)=\{\xi_2,\xi_3,...,\xi_r,x_1^+\}\subset End(T_{k+1}^\prime)$. Since $|\mu^+\overrightarrow{M_1}y_1^-|\le \lambda-1$, by Claim 10, $N(x_1^+)\subseteq V(H^\prime)$. Further, by Claim 16, $N(v)\subseteq V(H^\prime)$ for each $v\in End(H^\prime)$. Observing also that $|End(H^\prime)|=r\le m-1$, we can use the induction hypothesis, that is $|H^\prime|\ge \sigma_r+1$.

Next, if $d_{H}(x_1)=1$ and $w_i^\ast\not=x_1$, then $H^\prime=H\cup x_1\overrightarrow{M_1}\mu\xi_1-w_1w_1^\ast$ and $T_{k+1}^\prime=T_{k+1}+\xi_1\mu-w_1w_1^\ast$. Clearly, $End(H^\prime)=\{\xi_2,...,\xi_r,w_1\}\subseteq End(T_{k+1}^\prime)$. Since $|\mu^+\overrightarrow{M_1}y_1^-|\le \lambda-1$, by Claim 10, $N(x_1^+)\subseteq V(H^\prime)$, implying (by Claim 16) that $N(v)\subseteq V(H^\prime)$ for each $v\in End(H^\prime)$. Then we can argue as in Case $d_{H}(x_1)\ge 2$.

Finally, assume that $d_{H}(x_1)=1$ and $w_i^\ast=x_1$, implying that $r=1$ and $H=\xi_1\overrightarrow{Q_1}w_1w_1^\ast$. Define $H^\prime=\xi_1\overrightarrow{Q_1}w_1x_1\overrightarrow{M_1}\mu\xi_1$. Clearly, $|H^\prime|\ge d(\xi_1)+1\ge \sigma_1+1=\sigma_r+1$.

Now define $F^\prime$ as follows.

If $d_{F}(y_1)\ge 2$ then $F^\prime=F$. Clearly, $End(F^\prime)\subseteq End(T_{k+1})$. By Claim 17, $N(\xi_i)\subseteq V(F)$ $(i=r+1,r+2,...,m)$. By the induction hypothesis, $|F|\ge \sigma_{m-r}+1$.

Now let $d_{F}(y_1)=1$. If $r=m-1$, that is $F=\xi_m\overrightarrow{Q_m}w_mw_m^\ast$ (where $w_m^\ast=y_1$) then clearly, $|F|\ge \sigma_1+1=\sigma_{m-r}+1$. Otherwise ($r\le m-2$) $F$ has a branch vertex $z_1$. Let $\overrightarrow{R_1}=z_1\overrightarrow{R_1}y_1$ be the path connecting $z_1$ to $y_1$ in $F$. Choose $z_1$ so that $R_1$ is as short as possible. Put $F^{(1)}=F-V(R_1-z_1)$. If $N(\xi_i)\subseteq V(F^{(1)})$ for each $i\in\{r+1,r+2,...,m\}$ then by induction hypothesis, $|F|\ge|F^{(1)}|\ge \sigma_{m-r}+1$. Otherwise, let $\xi_jz_2\in E(G)$ for some $j\in\{r+1,...,m\}$, say $j=r+1$, and $z_2\in R_1-z_1$. Put $F^{(2)}=F^{(1)}\cup z_1\overrightarrow{R_1}z_2+\xi_{r+1}z_2-z_1z_1^+$. Clearly, $End(F^{(2)})=\{\xi_{r+2},...,\xi_m,z_1^+\}$. If $N(z_1^+)\subseteq V(F^{(2)})$ then by the induction hypothesis, $|F|\ge |F^{(2)}|\ge\sigma_{m-r}+1$. Otherwise let $z_1^+z_3\in E(G)$ for some $z_3\in z_2^+\overrightarrow{R_1}y_1$. Put
$$
F^{(3)}=F^{(2)}\cup z_2^+\overrightarrow{R_1}z_3+z_1^+z_3-z_2z_2^+.
$$
Define $F^{(1)}, F^{(2)},...$ and $z_1,z_2,...$ as long as possible and let $h$ be the maximum integer such that $N(z_{h-1}^+)\subseteq End(F^{(h)})$. It follows that $N(v)\subseteq V(F^{(h)})$ for each $v\in End(F^{(h)})\subset End(T_{k+1})$. By induction hypothesis, $|F|\ge|F^{(h)}|\ge\sigma_{m-r}+1$. So, in any case, $|H|\ge \sigma_r+1$ and $|F|\ge\sigma_{m-r}+1$. Then
$$
|T_m|\ge|H|+|F|\ge \sigma_r+\sigma_{m-r}+2\ge \sigma_m+1,
$$
and (2) holds by Claim 5.\\

\textbf{Case 2.2.3}. $N_{D_1}(\xi_i)=\emptyset$ $(i=1,...,m)$.

By the hypothesis (Case 2.2), $|D_1|\le\lambda-1$. By Claim 10, $N(\xi_i)\subseteq V(H)$ for each $i\le r$, and $N(\xi_i)\subseteq V(F)$ for each $i\ge r+1$. Then we can argue as in Case 2.2.2.

\noindent Institute for Informatics and Automation Problems\\ National Academy of Sciences\\
P. Sevak 1, Yerevan 0014, Armenia\\
E-mail: zhora@ipia.sci.am


\begin{thebibliography}{10}

\bibitem{[1]} J.A. Bondy and U.S.R. Murty, Graph Theory with Applications, Macmillan, London and Elsevier, New York (1976).

\bibitem{[2]} H. Broersma and H. Tuinstra, Independence trees and Hamilton cycles, J. Graph Theory 29 (1998) 227-237.

\bibitem{[3]} M. Las Vergnas, Sur une propri\'{e}t\'{e} des arbres maximaux dans un graphe, C.R. Acad.Sci.Paris S\'{e}r. A 272 (1971) 1297-1300.

\bibitem{[4]} Zh.G. Nikoghosyan, Two theorems on spanning trees, Uchenie Zapiski EGU (Scientific Transactions of the Yerevan StatUniversity), Ser. Matematika, issue 3 (1976) 3-6 (in Russian).

\bibitem{[5]} O. Ore, A note on hamiltonian circuits, Am. Math. Month. 67 (1960) 55.



\end{thebibliography}
\end{document}